\input amstex
\documentstyle{amsppt}
\NoBlackBoxes \magnification=\magstep1 \vsize=22 true cm \hsize=16.1
true cm \voffset=1 true cm
\pageno=0

\redefine\BMO{\operatorname{BMO}}
\redefine\VMO{\operatorname{VMO}}
\redefine\WP{\operatorname{WP}}
\redefine\Aff{\operatorname{Aff}}

\redefine\Diff{\operatorname{Diff}}

\topmatter

\title     Weil-Petersson Teichm\"{u}ller space
\endtitle
\title     Weil-Petersson Teichm\"{u}ller space III: dependence of Riemann mappings for Weil-Petersson curves
\endtitle
\author  Yuliang Shen  \quad Li Wu
\endauthor
\affil  Department of Mathematics, Soochow  University, Suzhou 215006, P.\,R. China
\endaffil

\abstract The classical Riemann mapping theorem implies that there exists a so-called Riemann mapping which takes the upper half plane onto the left domain bounded by a Jordan curve in the extended complex plane. The primary purpose of the paper is to study the basic problem: how does a Riemann mapping depend on the corresponding Jordan curve? We are mainly concerned with those Jordan curves in the  Weil-Petersson class, namely, the corresponding Riemann mappings can be quasiconformally extended to the whole plane with Beltrami coefficients being square integrable under the Poincar\'e metric. After giving a geometric characterization of a Weil-Petersson curve, we endow the space of all normalized Weil-Petersson curves with a new real Hilbert manifold structure
in  a geometric manner and show that this new structure  is topologically equivalent to the standard complex Hilbert manifold structure, which implies that an appropriately chosen Riemann mapping  depends   continuously on a Weil-Petersson curve (and vice versa). This can be considered as the first result about the continuous dependence of  Riemann mappings on non-smooth Jordan curves.

\vskip 0.2 cm


\vskip 0.2 cm


\endabstract

\thanks Research  supported by  the National Natural Science Foundation of China (Grant No. 11631010).
\endthanks

\thanks Email: ylshen\@suda.edu.cn\quad wuli187\@126.com
\endthanks
\endtopmatter

\document

\newpage

\head 1 Background \endhead

This is the third one of a series of papers (see [Sh], [ST]) which deal with  the Weil-Petersson geometry  theory of the universal Teichm\"uller space, a topic which  is important in Teichm\"uller theory (see [TT]) and has wide applications to various areas such as mathematical physics (see [BR1-2], [Ki], [KY], [RSW1-3]),  differential equation and   computer vision (see [FN], [Fi], [GMR], [GR], [Ku], [SM]) and the theory of Schramm-Loewner Evolutions (SLE) (see [VW], [Wa1-2]).

We begin with the well-known Riemann mapping theorem: Let $\Omega$ be the left domain bounded by a Jordan curve $\Gamma$ passing through the point at infinity in the extended complex plane $\hat\Bbb C$. Then there exists a univalent analytic function $f$ which maps the upper half plane $\Bbb U\doteq\{z=x+iy: y>0\}$ conformally onto $\Omega$ with $f(\infty)=\infty$. There also exists a univalent analytic function $g$ which maps the lower half plane $\Bbb U^*\doteq\{z=x+iy: y<0\}$ conformally onto the right domain $\Omega^*$ bounded by $\Gamma$ with $g(\infty)=\infty$. Both $f$ and $g$ are uniquely determined up to an affine mapping $z\mapsto az+b$ with $a>0$, $b\in \Bbb R$, the real line. $f$ and $g$ determine an increasing homeomorphism $h:\Bbb R\to\Bbb R$ by $h=f^{-1}\circ g$, which  is called a conformal sewing mapping of the curve $\Gamma$. $h$  is uniquely determined up to two affine mappings.  The primary purpose of this paper is to study the basic problem: how do the  mappings $f$, $g$ and $h$ depend on a  Jordan curve $\Gamma$? This problem was investigated in an important paper by Coifman-Meyer [CM] and is  a good example of a problem in nonlinear Fourier analysis,  as explained by Semmes [Se4].

Recall that a Jordan curve $\Gamma$ passing through the point at infinity  is a chord-arc (or Lavrentiev) curve with constant $k\ge 0$ if it is locally rectifiable and $$|s_1-s_2|\le (1+k)|z(s_1)-z(s_2)|\tag 1.1$$
for all $s_1\in \Bbb R$ and $s_2\in\Bbb R$, where $z(s)$ is a parametrization of $\Gamma$ by the arc-length $s\in\Bbb R$ (see [La], [Po2]). Coifman-Meyer [CM] showed that a Riemann mapping $f$ depends on $\Gamma$ real-analytically when $\Gamma$ is a chord-arc curve passing $\infty$ (see  [Wu] for an analogous result for bounded chord-arc curves).  To make this precise, let $\Gamma$ be a chord-arc curve passing through $0$ and $\infty$, and $z(s)$ be the (unique) arc-length parametrization of $\Gamma$ with $z(0)=0$. David [Da] showed that  there exists some function $b$ in $\BMO_{\Bbb R}$ (or more precisely, $\BMO_{\Bbb R}/\Bbb R$), the space of all real-valued functions of bounded mean oscillation on the real line (see [FS], [Gar], [Po2],
[Zh] and section 3 below),   such that $z'(s)=e^{ib(s)}$, and these BMO functions $b's$ form an open subset $\Cal L$ of $\BMO_{\Bbb R}/\Bbb R$. A Riemann mapping $f:\Bbb U\to\Omega$ induces an increasing homeomorphism $h_1:\Bbb R\to\Bbb R$ by $f\circ h_1=z$. A classical result of Lavrentiev [La] implies that $h_1$ is locally absolutely continuous so that $h'_1$ belongs to the class of weights $A^{\infty}$ introduced
by Muckenhoupt (see [CF], [Gar]), in particular, $\log  {h'_1}$ is a BMO function. The precise statement of the result of Coifman-Meyer [CM] is: The correspondence $b\mapsto\log h'_1$ induces a well-defined  real-analytic map from $\Cal L$ into $\BMO_{\Bbb R}/\Bbb R$. A different approach to this result was given later by Semmes [Se3] (see also [Se1]).

To see how a Riemann mapping $f$ itself, not just the induced mapping $h_1$ by $f\circ h_1=z$, depends on the curve $\Gamma$ when it is a chord-arc curve, we recall a result of Pommerenke [Po1] (see also [Zi]), which says that $\log f'$ belongs to {\rm BMOA}, the space of
analytic functions in $\Bbb U$ of bounded mean oscillation (see [FS], [Gar]). From
$$f\circ h_1=z\Rightarrow (f'\circ h_1)h'_1=z'=e^{ib}\Rightarrow \log (f'\circ h_1)+\log h'_1=ib\Rightarrow\log f'=(ib-\log h'_1)\circ h_1^{-1},$$
it is not clear how the mapping $f$ depends on the curve $\Gamma$ (or on the function $b$), although $\log h'_1$ depends real analytically on $b$. For example, we do not know whether $\log f'$ depends continuously on $b$. It is also not clear how a conformal sewing mapping $h$ depends on the curve $\Gamma$ (or on the function $b$) when it is a chord-arc curve. Actually,  Katznelson-Nag-Sullivan
[KNS] asked whether $\log (h^{-1})'$ depends continuously on $b$ for a chord-arc curve $\Gamma$. Anyhow, by means of some results in our paper [SWe] (see also [AZ]), we conclude  that, under some  normalized conditions, $\log f'$ depends continuously on $\log(h^{-1})'$ and $\log(h^{-1})'$ depends continuously on $\log f'$, which implies that  the continuous dependence of $\log f'$ on $\Gamma$ (or $b$) would imply  the continuous dependence of $\log(h^{-1})'$ on $\Gamma$ (or $b$), and vice versa{\footnote{We conjecture that in general neither $\log f'$ nor $\log(h^{-1})'$ depends continuously on $\Gamma$ (or $b$).}}.

Since it is not clear  how a Riemann mapping $f$ (or $g$) depends on the curve $\Gamma$ even when it is a chord-arc curve, it is desirable to find a subclass of chord-arc curves on which a Riemann mapping does depend continuously. In this paper, we will consider an important sub-class of chord-arc curves, which we call them  Weil-Petersson curves. We say a Jordan curve $\Gamma$ passing through $\infty$ is  a Weil-Petersson curve if a Riemann mapping $g$,  which maps the lower half plane $\Bbb U^*$ conformally onto the right domain $\Omega^*$ bounded by $\Gamma$ with $g(\infty)=\infty$, has  a quasiconformal extension to the whole plane whose Beltrami coefficient is square integrable in the Poincar\'e metric. This class of Jordan curves and its Teichm\"uller space has been much investigated in recent years (see [Cu], [Fi], [GGPPR], [GR], [Sh], [ST], [STW], [TT] and section 2 below).
However, it is still an open problem how to give a geometric characterization of a Weil-Petersson curve without using a Riemann mapping $g$ or its quasiconformal extensions. An analogous problem was proposed by Takhtajan-Teo [TT] for bounded Weil-Petersson curves and a partial answer was given by Gallardo-Guti\'errez,  Gonz\'alez,  P\'erez-Gonz\'alez,  Pommerenke and  R\"atty\"a [GGPPR] in this case.

As a preliminary result we will first show that for the  arc-length parametrization $z(s)$, $z(0)=0$, of a Weil-Petersson curve $\Gamma$ passing through the points $0$ and $\infty$,   there exists some function $b$ in the real Sobolev class $H^{\frac 12}_{\Bbb R}/\Bbb R$  such that $z'(s)=e^{ib(s)}$, and these $H^{\frac 12}$ functions $b's$ form an open subset  of $H^{\frac 12}_{\Bbb R}/\Bbb R$. Then we will come back to our main goal of this paper and show that for those Weil-Petersson curves $\Gamma$ with $z(0)=0$, $z(1)>0$, $z(\infty)=\infty$, an appropriately chosen Riemann mapping $f$ (or $g$) and the corresponding conformal sewing mapping $h$ depend   continuously on  $\Gamma$ (and vice versa). The precise statements of these results will be given in section 2. As far as we know, this is the first result on continuous dependence of the Riemann mappings and conformal sewing mappings on non-smooth Jordan curves.

In the paper,  $C$,
$C_1$, $C_2$ $\cdots$ will denote universal constants that might
change from one line to another, while $C(\cdot)$, $C_1(\cdot)$,
$C_2(\cdot)$ $\cdots$ will denote constants that depend only on the
elements put in the brackets.
 The notation $A\lesssim B$ $(A\gtrsim B)$ means that there is a positive constant $C$ independent of $A$ and $B$ such that
$A\le CB$ $ (A\ge CB)$. The notation $A \asymp B$ means both $A\lesssim B$ and $A\gtrsim B$.

 \head 2 Introduction and statement of results
\endhead

In this section, we will give some basic definitions and results on the universal Teichm\"uller space and the Weil-Petersson Teichm\"uller space (see the books [Ah], [GL], [Le], [Na] and the papers [Sh], [ST], [TT] for more details). We will also state the main results of this paper.

\vskip 0.2 cm
\noindent {\bf 2.1 Various models of the universal Teichm\"uller space} \quad
Let  $M(\Bbb U)$ denote  the open
unit ball of the Banach space $L{}^{\infty}(\Bbb U)$ of
essentially bounded measurable functions on the upper half plane $\Bbb U$ in the complex plane $\Bbb C$. For $\mu\in
M(\Bbb U)$, let $f^{\mu}$ be the unique  quasiconformal mapping from $\Bbb U$ onto itself which has complex dilatation $\mu$ and fixes the points $0$, $1$ and $\infty$, and  $f_{\mu}$ be the unique  quasiconformal mapping on the
extended plane $\hat {\Bbb C}$ which has  complex dilatation
$\mu$ in $\Bbb U$,  is conformal in the lower half plane $\Bbb U^*$ and fixes the points $0$, $1$ and $\infty$. We say two elements $\mu$
and $\nu$ in $M(\Bbb U)$ are equivalent, denoted by $\mu\sim\nu$,
if $f^{\mu}|_{\Bbb R}=f^{\nu}|_{\Bbb R}$, or equivalently, $f_{\mu}|_{\Bbb U^*}=f_{\nu}|_{\Bbb U^*}$. Then
$T=M(\Bbb U)/_{\sim}$ is the Bers model of the universal
Teichm\"uller space. We let $\Phi$ denote the natural  projection
 from $M(\Bbb U)$ onto $T$ so that $\Phi(\mu)$ is the equivalence class $[\mu]$. $[0]$ is called the base point of $T$.

An  increasing     homeomorphism $h$ from  the real line $\Bbb R$  onto itself is said to be quasisymmetric if there exists a (least) positive constant $C(h)$, called the quasisymmetric constant of $h$,  such
that $|h(I_1)|\le C(h)|h(I_2)|$
for all pairs of adjacent intervals  $I_1$ and $I_2$ on $\Bbb R$ with the same length $|I_1|=|I_2|$. Beurling-Ahlfors [BA] proved that an increasing homeomorphism $h$ from  the real line $\Bbb R$  onto itself is quasisymmetric if and only if there
exists some quasiconformal homeomorphism  of $\Bbb U$ onto itself
which has boundary values $h$. A Jordan curve $\Gamma$ passing through $\infty$ is said to be a quasicircle if it is the image of the extended real line $\hat\Bbb R=\Bbb R\cup\{\infty\}$ under a quasiconformal mapping on the whole plane. It is easy to see that a Jordan curve $\Gamma$ passing through $\infty$ is  a quasicircle if and only if a Riemann mapping $f$ (or $g$) can be extended to a quasiconformal mapping on the whole plane. Therefore, the universal Teichm\"uller space $T$ can also be defined as:

\noindent
$\bullet$ The set of all  quasisymmetric homeomorphisms of the real line onto itself with $0$, $1$ and $\infty$ fixed ($[\mu]\mapsto f^{\mu}|_{\Bbb R}$).

\noindent
$\bullet$ The set of all  conformal mappings on the lower half plane $\Bbb U^*$ which can be quasiconformally extended to  the whole plane with the points $0$, $1$ and $\infty$ fixed ($[\mu]\mapsto f_{\mu}|_{\Bbb U^*}$).

\noindent
 $\bullet$  The set of all  quasicircles  through the points $0$, $1$ and $\infty$ ($[\mu]\mapsto f_{\mu}(\hat\Bbb R)$).

It is known that the universal Teichm\"uller space $T$ is an infinite dimensional
complex Banach manifold. To make this precise, we first recall some important Banach spaces.   Let  $D$ be an
arbitrary simply connected domain in the extended complex plane
$\hat\Bbb C$ which is hyperbolic in the sense that it is conformally equivalent to the upper half plane. The hyperbolic metric $\lambda_{D}(z)|dz|$ (with curvature constantly equal to $-1$) in $D$ is
defined by
$$\lambda_{D}(f(z))|f'(z)|=\frac{1}{y},\quad z=x+iy\in\Bbb U,\tag 2.1$$
where  $f: \Bbb U\rightarrow D$ is any conformal
mapping.
Let $B_2(D)$ denote the Banach space of functions $\phi$
holomorphic in $D$ with norm
$$
\|\phi\|_{B_2(D)}\doteq\sup_{z\in D}|\phi(z)|\lambda^{-2}_{D}(z), \tag 2.2
$$
and  $ \Cal B(D)$ the Banach space of functions $\phi$
holomorphic in $D$ with finite norm
$$\|\phi\|_{ \Cal B(D)}\doteq\left(\frac{1}{\pi}\iint_{D}|\phi(z)|^2\lambda^{-2}_{D}dxdy\right)^{\frac 12}.\tag 2.3$$
Then,  $ \Cal B(D)\subset B_2(D)$, and the
inclusion map is continuous (see [Zh]).

Now  we consider the map $S: M(\Bbb U)\to B_2(\Bbb U^*)$ which sends $\mu$ to
 the Schwarzian derivative of $f_{\mu}|_{\Bbb U^*}$. Recall that for any locally univalent function $f$, its Schwarzian derivative $S_f$ is defined by
 $$S_f\doteq N'_f-\frac 12N^2_f,\quad N_f\doteq(\log f')'.\tag 2.4$$
 $S$ is a holomorphic split submersion onto its image, which
descends down to a map $\beta: T\to B_2(\Bbb U^*)$  known as the Bers
embedding. Via the Bers embedding, $T$ carries a natural complex Banach manifold
structure so that  $\Phi$  is a holomorphic split submersion.

Besides the Schwarzian derivative model, the universal Teichm\"uller space has another important model, the pre-logarithmic derivative model.  In the unit disk case, the  pre-logarithmic derivative model of the universal  Teichm\"uller spaces was much investigated  (see [AGe], [Po2], [Zhu]). Here we consider the upper half plane case. Let
 $B_1(D)$
denote the Bloch space of functions $\phi$ holomorphic in a hyperbolic simply connected domain $D$
with semi-norm
$$
\|\phi\|_{B_1(D)}\doteq\sup_{z\in D}|\phi'(z)|\lambda^{-1}_{D}(z), \tag 2.5
$$
and
 $\Cal D(D)$
denote the Dirichlet  space of functions $\phi$ holomorphic in $D$
with semi-norm
$$
\|\phi\|_{\Cal D(D)}\doteq\left(\frac{1}{\pi}\iint_{D}|\phi'(z)|^2dxdy\right)^{\frac 12}. \tag 2.6
$$
It is known that   $ \Cal D(D)\subset B_1(D)$, and the
inclusion map is continuous (see [Zh]). It is also known that, for each holomorphic function $\phi$ on $D$, $\phi''\in \Cal B(D)$ if $\phi\in \Cal D(D)$, and $\phi''\in B_2(D)$ if $\phi\in  B_1(D)$. The converse is also true,  with some normalized conditions of $\phi$ at $\infty$ whenever $D$ is not a bounded  domain (see [ST], [STW]).

Now Koebe distortion theorem implies that $\log (f_{\mu}|_{\Bbb U^*})'\in B_1(\Bbb U^*)$ for $\mu\in M(\Bbb U)$. Furthermore, the map $L$ induced by the correspondence $\mu\mapsto \log (f_{\mu}|_{\Bbb U^*})'$ is  a continuous map from $M(\Bbb U)$ into $B_1(\Bbb U^*)$ (see [Le]). Actually, $L: M(\Bbb U)\to B_1(\Bbb U^*)$ is even holomorphic (see [Ha]).

\vskip 0.2 cm

\noindent {\bf 2.2 Various models of the Weil-Petersson Teichm\"uller space} \quad
Now we define
 the Weil-Petersson Teichm\"uller space. We denote by
$\Cal L^{\infty}(D)$ the Banach space of all essentially bounded
measurable functions $\mu$ on a simply connected domain $D$  with norm
$$\|\mu\|_{\WP}\doteq\|\mu\|_{\infty}+\left(\frac{1}{\pi}\iint_{D}|\mu(z)|^2\lambda^2_{D}(z)dxdy\right)^{\frac 12}.\tag 2.7$$
 Set $\Cal M(\Bbb U)=M(\Bbb U)\cap \Cal L^{\infty}(\Bbb U)$. Then $\Cal T=\Cal M(\Bbb U)/_{\sim}$ is known as  the Weil-Petersson Teichm\"uller space. Actually, $ \Cal T$ is the base point component of the universal Teichm\"uller space under the complex Hilbert manifold structure introduced by  Takhtajan-Teo [TT]. Under the Bers projection $S: M(\Bbb U)\to
B_2(\Bbb U^*)$, $S(\Cal M(\Bbb U))=S(M(\Bbb U))\cap  \Cal B(\Bbb U^*)$ (see [Cu], [TT]). More precisely, $S:\Cal M(\Bbb U)\to \Cal B(\Bbb U^*) $ is a holomorphic split submersion onto its image, which
induces  a natural complex Hilbert  manifold on $\Cal T$  so that  $\Phi: \Cal M(\Bbb U)\to\Cal T$  is a holomorphic split submersion.
Very recently, we proved that under the pre-logarithmic derivative projection $L: M(\Bbb U)\to
B_1(\Bbb U^*)$, $L(\Cal M(\Bbb U))=L(M(\Bbb U))\cap  \Cal D(\Bbb U^*)$, and $L: \Cal M(\Bbb U)\to \Cal D(\Bbb U^*)$ is holomorphic (see [STW]).

We proceed to  introduce the quasisymmetric homeomorphism model of the Weil-Petersson Teichm\"uller space. For simplicity, we say a quasiconformal mapping on a hyperbolic simply connected domain $D$ is of the Weil-Petersson class if its Beltrami coefficient is in $\Cal L^{\infty}(D)$. A quasiconformal mapping $f$ on the whole plane is said to belong to the Weil-Petersson class (with respect to the real line) if both $f|_{\Bbb U}$ and $f|_{\Bbb U^*}$ are of the  Weil-Petersson class. A sense preserving  homeomorphism $h$ of the  real line $\Bbb R$ onto itself is said to belong the Weil-Petersson class if it can be extended a Weil-Petersson quasiconformal mapping to  the upper half plane $\Bbb U$. We denote by $\WP(\Bbb R)$ the class of all Weil-Petersson homeomorphisms on $\Bbb R$, and by $\WP_0(\Bbb R)$ be those with the points $0$, $1$ and $\infty$ fixed. Then we have the following result.

 \proclaim{Proposition 2.1 ([ST], [STW])} An increasing homeomorphism $h$ from the real line $\Bbb R$ onto itself belongs to the Weil-Petersson class if and only if $h$ is locally absolutely continuous with $\log h'\in H^{\frac 12}$. Moreover, when $\log h'\in H^{\frac 12}$, $h$ can be can be extended to a Weil-Petersson  quasiconformal mapping  of the upper half plane $\Bbb U$ onto itself which is also  bi-Lipschitz under the Poincar\'e metric $\lambda_{\Bbb U}(z)|dz|$.
 \endproclaim

  Recall that the Sobolev  class $H^{\frac 12}$ ($H^{\frac 12}_{\Bbb R}$) on the real line $\Bbb R$  is the collection of all locally integrable (real) functions $u$ with
$$\|u\|^2_{H^{\frac 12} }\doteq\frac{1}{4\pi^2}\int_{\Bbb R}\int_{\Bbb R}\frac{|u(s)-u(t)|^2}{(s-t)^2}dsdt<+\infty.\tag 2.8$$
Proposition 2.1 was proved in our papers [ST], [STW]. An analogous result on the unit circle was  proved earlier by the first author [Sh] (see also [WHS]), which solves a problem proposed by Takhtajan-Teo in 2006 (see page 68 in [TT] and also [Fi], [GR]). The following result (see [ST], [STW]) says that the  normalized Weil-Petersson class $\WP_0(\Bbb R)$,   the quasisymmetric homeomorphism model of the  Weil-Petersson Teichm\"uller space $\Cal T$, can be endowed with a real Hilbert manifold structure from $H_{\Bbb R}^{\frac 12}/\Bbb R$   by the correspondence $h\mapsto\log h'$, which is real analytically equivalent to the standard complex Hilbert manifold structure on $\Cal T$   given by Takhtajan-Teo [TT]. Both Propositions 2.1 and 2.2 will play an important role in our later discussion.

\proclaim{Proposition 2.2 ([ST])} The correspondence $h\mapsto\log h'$  induces a real analytic map $\Psi$  from the  normalized Weil-Petersson class $\WP_0(\Bbb R)(=\Cal T)$ onto the real Sobolev space $H_{\Bbb R}^{\frac 12}/\Bbb R$  whose inverse $\Psi^{-1}$ is also real analytic.

\endproclaim

Before we state the main results, we summarize that the Weil-Petersson Teichm\"uller space $\Cal T$ can also be defined in the following ways:

\noindent
$\bullet$ The set $\WP_0(\Bbb R)$ of all normalized Weil-Petersson  homeomorphisms with $0$, $1$, $\infty$ fixed ($[\mu]\mapsto f^{\mu}|_{\Bbb R}$).

\noindent
$\bullet$ The set of all conformal mappings $g$ on the upper half plane $\Bbb U^*$ which can be extended to Weil-Petersson quasiconformal mappings on the whole plane and satisfies the following normalized conditions ($[\mu]\mapsto g_{\mu}|_{\Bbb U^*}$):
 $$g(0)=0,\, g(\infty)=\infty,\, g(1)>0,\, \int_0^1|g'(t)|dt=1.\tag 2.9 $$

\noindent
 $\bullet$  The set of all normalized Weil-Petersson curves on the whole plane ($[\mu]\mapsto g_{\mu}(\hat\Bbb R)$).

\noindent For later purposes, here we have used some  normalized conditions different from the universal Teichm\"uller space case. A Weil-Petersson curve $\Gamma$ is called normalized if it passes through $0$ and $\infty$ and the arclength parametrization  $z=z(s)$   of  $\Gamma$ with  $z(0)=0$ satisfies $z(1)>0$. For $\mu\in
\Cal M(\Bbb U)$,   $g_{\mu}$ is the unique  quasiconformal mapping on the
extended plane $\hat {\Bbb C}$ which has  complex dilatation
$\mu$ in $\Bbb U$,  is conformal in  $\Bbb U^*$ and with the normalized conditions (2.9). It is easy to see that  $\mu$
and $\nu$ in $\Cal M(\Bbb U)$ are equivalent if and only if  $g_{\mu}|_{\Bbb U^*}=g_{\nu}|_{\Bbb U^*}$. Actually, it holds that $g_{\mu}=g_{\mu}(1)f_{\mu}$.

\vskip 0.2 cm

\noindent {\bf 2.3 Statement of main results} \quad As stated in section 1, an open problem   is to give a geometric characterization of a Weil-Petersson curve without using a Riemann mapping $f$ (or $g$) or its quasiconformal extensions. A basic geometric notion to a locally rectifiable curve is an  arc-length parametrization. Therefore, a natural question is to characterize an arclength parametrization of a Weil-Petersson curve. As a preliminary result of the paper  we will first show

\proclaim{Theorem 2.1} Let $\Gamma$ be a normalized Weil-Petersson curve  and $z=z(s)$ be the  arc-length parametrization  of  $\Gamma$ with  $z(0)=0$. Then there exists some function $b$ in the real Sobolev class $H^{\frac 12}_{\Bbb R}/\Bbb R$  such that $z'(s)=e^{ib(s)}$. Moreover, the set  $\hat \Cal T$ of these $H^{\frac 12}$ functions $b's$ is an open subset of $H^{\frac 12}_{\Bbb R}/\Bbb R$.
\endproclaim

Actually, we have the following geometric characterization of a Weil-Petersson curve by means of the arc-length parametrization under the geometric assumption of chord-arc property. A similar result also holds   for bounded Weil-Petersson curves. In fact, Bishop [Bi] has obtained  several geometric characterizations of a bounded Weil-Petersson curve very recently{\footnote{After an earlier version of this manuscript was posted on arXiv [SWu], we learned from Tim Mesikepp that Bishop [Bi] obtained  various geometric characterizations of a bounded Weil-Petersson curve. The
authors would like to thank Tim Mesikepp for calling this reference to their
attention.}}.

\proclaim{Theorem 2.2} Let $\Gamma$ be a locally rectifiable Jordan curve passing through $\infty$ and  $z=z(s)$ be an  arc-length parametrization  of  $\Gamma$. Then $\Gamma$ is a Weil-Petersson curve if and only if $\Gamma$ is  a  chord-arc curve  and there exists some function $b$ in the real Sobolev class $H^{\frac 12}_{\Bbb R}/\Bbb R$  such that $z'(s)=e^{ib(s)}$. In other words, $\hat \Cal T=\Cal L\cap H^{\frac 12}_{\Bbb R}/\Bbb R$.
\endproclaim

Theorem 2.1 implies that the set of all normalized Weil-Petersson curves, a model of the Weil-Petersson Teichm\"uller space $\Cal T$, can be  endowed with  a real Hilbert manifold structure in a geometric manner by the correspondence $\Gamma\mapsto b$.  We will show that this new real Hilbert manifold structure is topologically equivalent to the  standard complex Hilbert manifold structure  given by Takhtajan-Teo [TT]. To be precise, we introduce some notations. For a normalized Weil-Petersson curve $\Gamma$ with arc-length parametrization $z(s)$, $z(0)=0$, we  denote  by $f_{\Gamma}$  the unique Riemann mapping which takes $\Bbb U$ onto the left domain $\Omega$ bounded  $\Gamma$ with the normalized conditions (2.9), that is, $f_{\Gamma}(s)=z(s)$ for $s=0, 1, \infty$. Similarly, we denote by $g_{\Gamma}$  the unique Riemann mapping which takes $\Bbb U^*$ onto the right domain $\Omega^*$ bounded  $\Gamma$ such that $g_{\Gamma}(s)=z(s)$ for $s=0, 1, \infty$.  Denote by $h_{\Gamma}$ the unique conformal sewing mapping determined by $f_{\Gamma}$ and $g_{\Gamma}$, that is, $h_{\Gamma}=f_{\Gamma}^{-1}\circ g_{\Gamma}$. For $b\in\hat\Cal T$, we may assume without loss of generality that $\int_0^1e^{ib(t)}dt>0$, and then  denote by $\Gamma_b$ the unique normalized Weil-Petersson curve whose arc-length parametrization $z_b$ with $z_b(0)=0$ satisfies $z'_b=e^{ib}$, namely,
$$z_b(s)=\int_0^se^{ib(t)}dt.\tag 2.10$$
Finally, for $b\in\hat\Cal T$, we set $f_b=f_{\Gamma_b}$, $g_b=g_{\Gamma_b}$, and $h_b=h_{\Gamma_b}$.

Now we can state the main results  of this paper. They say as we have promised that an appropriately chosen Riemann mapping $f$ (or $g$) and the corresponding conformal sewing mapping $h$ depend   continuously on a normalized  Weil-Petersson curve, and vice versa.

\proclaim{Theorem 2.3} The correspondence $b\mapsto h_b$ induces a homeomorphism from $\hat\Cal T$ onto the normalized Weil-Petersson class $\WP_0(\Bbb R)\,(=\Cal T)$, or equivalently, the correspondence $b\mapsto \log h'_b$ induces a homeomorphism from $\hat\Cal T$ onto the real Sobolev space $H_{\Bbb R}^{\frac 12}/\Bbb R$.
\endproclaim

Theorem 2.3 can also be restated as

\proclaim{Theorem 2.4} The correspondence  $b\mapsto g_b$ induces a homeomorphism from $\hat\Cal T$ onto (the conformal mapping model of) the Weil-Petersson Teichm\"uller space $\Cal T$. Equivalently, the correspondence $b\mapsto \log g'_b$ induces a homeomorphism from $\hat\Cal T$ onto its image in $\Cal D(\Bbb U^*)/\Bbb C$, or equivalently, the correspondence $b\mapsto \log f'_b$ induces a homeomorphism from $\hat\Cal T$ onto its image in $\Cal D(\Bbb U)/\Bbb C$.
\endproclaim

Since the Weil-Petersson Teichm\"uller space $\Cal T$ is contractible, we obtain

\proclaim{Corollary 2.1} The arc-length parametrization space $\hat\Cal T$ of the normalized Weil-Petersson curves is contractible.
\endproclaim

\noindent
{\bf Remark.}\quad  Recall that $\hat \Cal T=\Cal L\cap H^{\frac 12}_{\Bbb R}/\Bbb R$. It is not known whether the arc-length parametrization space $\Cal L$ of the normalized chord-arc curves is contractible. Actually, it is even not known whether $\Cal L$ is connected. This is known to be a difficult open problem (see [AGo], [AZ], [CM], [Se1]).

During the proof of Theorems 2.3 and 2.4, we will obtain a much stronger result than Corollary 2.1.  Let $\hat\Cal T_e$ denote the set of all elements $u$ in $H^{\frac 12}$ (or more precisely, $H^{\frac 12}/\Bbb C$) such that the function $\gamma_u$ defined by
$$\gamma_u(x)=\int_0^xe^{iu(t)}dt,\quad x\in\Bbb R,\tag 2.11$$
is a homeomorphism from the extended real line $\hat\Bbb R$ onto a Weil-Petersson curve. Clearly, the arc-length parametrization space $\hat\Cal T$ of the normalized Weil-Petersson curves comprises precisely the real-valued functions in $\hat\Cal T_e$, that is, $\hat\Cal T=\hat\Cal T_e\cap H_{\Bbb R}^{\frac 12}/\Bbb R$.

\proclaim {Theorem 2.5} The parametrization space $\hat\Cal T_e$ is a contractible domain in $H^{\frac 12}/\Bbb C$.
\endproclaim

\noindent
{\bf Remark.}\quad Proposition 2.1 says that an increasing homeomorphism $h$ from the real line $\Bbb R$ onto itself can be extended a Weil-Petersson quasiconformal mapping to the whole plane  if and only if $h$ is locally absolutely continuous with $\log h'\in H^{\frac 12}$. We will generalize this result and show in section 5 (see Propositions 5.1-2 below) that a sense-preserving  homeomorphism $h$ on the real line $\Bbb R$  with $h(\infty)=\infty$ can be extended a Weil-Petersson quasiconformal mapping to the whole plane  if and only if $h$ is locally absolutely continuous such that  $\log h'\in H^{\frac 12}$ and it maps $\hat\Bbb R$ onto a Weil-Petersson curve containing $\infty$, or equivalently, $h=\gamma_u$ for some $u\in \hat\Cal T_e$ up to  an affine map.

 \head 3 BMO functions  revisited
\endhead

In order to prove Theorems 2.3-4, we need a construction concerning quasiconformal extensions of strongly quasisymmetric homeomorphisms introduced by Semmes [Se1-2], which relies heavily on BMO estimates. In this section we  recall  some basic definitions and results on BMO functions (see [Gar]).

A locally  integrable function $u\in L^1_{loc}(\Bbb R)$ is said to have bounded mean oscillation and belongs to the space $\BMO$ if
$$\|u\|_{\BMO}\doteq\sup\frac{1}{|I|}\int_I|u(t)-u_I|dt<+\infty,\tag 3.1$$
where the supremum is taken over all finite sub-intervals $I$ of $\Bbb R$, while $u_I$ is the average of $u$ on the interval $I$, namely,
$$u_I=\frac{1}{|I|}\int_Iu(t)dt.\tag 3.2$$ If $u$ also satisfies the condition
$$\lim_{|I|\to 0}\frac{1}{|I|}\int_I|u(t)-u_I|dt=0,$$
we say $u$ has vanishing mean oscillation and belongs to the space $\VMO$.  In the following, we denote by $\BMO_{\Bbb R}$ the set of all real-valued BMO functions.
 It is well known that  $H^{\frac 12}\subset \VMO$, and the inclusion map is continuous (see [Zh]). In fact,
$$
\align
\frac{1}{|I|}\int_I|u(t)-u_I|dt&=\frac{1}{|I|}\int_I\left|u(t)-\frac{1}{|I|}\int_Iu(s)ds\right|dt\\
&\le\frac{1}{|I|^2}\int_I\int_I|u(s)-u(t)|dsdt\\
&\le\frac{1}{|I|^2}\left(\int_I\int_I\frac{|u(s)-u(t)|^2}{(s-t)^2}dsdt\right)^{\frac 12}\left(\int_I\int_I(s-t)^2dsdt\right)^{\frac 12}\\
&\le \left(\int_I\int_I\frac{|u(s)-u(t)|^2}{(s-t)^2}dsdt\right)^{\frac 12},
\endalign
$$
which implies that $H^{\frac 12}\subset \VMO$, and $\|u\|_{\BMO}\lesssim \|u\|_{H^{\frac 12} }$.

We need some basic results on BMO functions. By the well-known theorem of John-Nirenberg for BMO functions (see [Gar]), there exist two universal positive constants $C_1$ and $C_2$ such that for any BMO function $u$, any subinterval $I$ of $\Bbb R$  and any $\lambda>0$, it holds that
 $$\frac{\left|\{t\in I:|u(t)-u_I|\ge\lambda\}\right|}{|I|}\le C_1\exp\left(\frac{-C_2\lambda}{\|u\|_{\BMO}}\right).\tag 3.3$$
  By  Chebychev's inequality, we obtain that for $u$ with $\|u\|_{\BMO}<C_2$,
 $$
 \align
 \frac{1}{|I|}\int_{I}(e^{|u-u_{I}|}-1)dt&=\frac{1}{|I|}\int_0^{\infty}\left|\{t\in I:|u-u_{I}|\ge\lambda\}\right|d(e^{\lambda}-1) \\
 &\le C_1\int_0^{\infty}e^\lambda\exp\left(\frac{-C_2\lambda}{\|u\|_{\BMO}}\right)d\lambda\tag 3.4\\
 &\le  \frac{C_1\|u\|_{\BMO}}{C_2-\|u\|_{\BMO}}.
 \endalign
 $$
Similarly, for any $p\ge 1$ we have
$$
  \frac{1}{|I|}\int_{I}|u-u_{I}|^pdt\lesssim C(p)\|u\|^p_{\BMO}. \tag 3.5
 $$

 \proclaim {Lemma 3.1} Let $\phi$ be a $C^{\infty}$  function on the real line which is supported on $[-1, 1]$ and satisfies $\int_{\Bbb R}\phi(x)dx=1$. Set $\phi_y(x)=|y|^{-1}\phi(|y|^{-1}x)$ for $y\neq 0$, and consider the convolution
 $$\phi_y\ast w(x)=\int_{\Bbb R}\phi_y(x-t)w(t)dt.\tag 3.6$$
 Suppose $v\in L^{\infty}(\Bbb R)$ and $|\phi_y\ast v|\ge \epsilon_0$ for some  $\epsilon_0>0$. Then for
 $$R_y(u)(x)=\frac{\phi_y\ast (vu)(x)}{\phi_y\ast v(x)}$$ it holds that $$|R_y(e^u)|\asymp |e^{R_y(u)}|\tag 3.7$$
 when $\|u\|_{\BMO}$ is small.
 \endproclaim

 \demo{Proof} Lemma 3.1 appeared implicitly in [Se3] though not stated in this form. For the convenience of later use, we write down  the detailed proof here (see [Se3], [ST]). Actually, not only  Lemma 3.1 itself, but also both of the estimates (3.8) and (3.9) below will be frequently used in section 6.

 For $x\in\Bbb R$ and $y>0$, consider $I=[x-y, x+y]$ so that
 $$u_I=\frac{1}{2y}\int^{x+y}_{x-y}u(t)dt.$$
Since $\int_{\Bbb R}\phi(x)dx=1$, which implies that $\int_{\Bbb R}\phi_y(x)dx=1$, and $R_y(1)\equiv 1$, we obtain
 $$
 \align
 |R_y(u)(x)-u_I|&=|R_y(u-u_I)(x)|\le\frac{1}{\epsilon_0}|\phi_y\ast (v(u-u_I))(x)|\\
 &\le \frac{C(\phi)\|v\|_{\infty}}{\epsilon_0} \frac{1}{|I|}\int_I|u(t)-u_I|dt\lesssim\|u\|_{\BMO}.\tag 3.8
 \endalign$$
 Now for any complex number $z$, it holds that $|e^z-1|\le e^{|z|}-1\le |z|e^{|z|}$, which goes as follows,
 $$|e^z-1|=\left|\sum_{n=1}^{+\infty}\frac{z^n}{n!}\right|\le \sum_{n=1}^{+\infty}\frac{|z|^n}{n!}=e^{|z|}-1,$$
 $$e^{|z|}-1=\sum_{n=1}^{+\infty}\frac{|z|^n}{n!}=|z|\sum_{n=1}^{+\infty}\frac{|z|^{n-1}}{n!}\le |z|\sum_{n=1}^{+\infty}\frac{|z|^{n-1}}{(n-1)!}=|z|e^{|z|}.$$
Then  we have
 $$
 \align
 &\frac{1}{|I|}\int_I|e^{u(t)-R_y(u)(x)}-1|dt\\
 &\le\frac{1}{|I|}\int_I|e^{u(t)-R_y(u)(x)}||u(t)-R_y(u)(x)|dt\\
 &\le\frac{|e^{u_I-R_y(u)(x)}|}{|I|}\int_I|e^{u(t)-u_I}|(|u(t)-u_I|+|u_I-R_y(u)(x)|dt.\\
 \endalign
 $$
 Using H\"older inequality, we conclude from (3.4), (3.5) and (3.8) that
 $$
 \frac{1}{|I|}\int_I|e^{u(t)-R_y(u)(x)}-1|dt\lesssim\|u\|_{\BMO}\tag 3.9$$
 when $\|u\|_{\BMO}$ is small. Noting that
 $$R_y(e^u)(x)-e^{R_y(u)(x)}=e^{R_y(u)(x)}R_y(e^{u-R_y(u)(x)}-1)(x),$$
 we obtain
 $$\align
 |R_y(e^u)(x)-e^{R_y(u)(x)}|&=|e^{R_y(u)(x)}||R_y(e^{u-R_y(u)(x)}-1)(x)|\\
 &\lesssim \frac{|e^{R_y(u)(x)}|}{|I|}\int_I|e^{u(t)-R_y(u)(x)}-1|dt,
 \endalign$$
 which implies by (3.9)  the required relation (3.7). \quad$\square$
 \enddemo

 \head 4 Proof of Theorems 2.1 and 2.2
\endhead
In this section, we will give  simple proof of Theorems 2.1 and 2.2. However, to prove Theorems 2.3 and 2.4, we need a concrete approach to the openness of $\hat\Cal T$, which will be given in section 6.

 Let $\Gamma$ be a locally rectifiable Jordan curve passing through $\infty$ and  $z=z(s)$ be an  arc-length parametrization  of  $\Gamma$. Let  $f$ map the upper half plane $\Bbb U$ conformally onto  the left domain $\Omega$ bounded by  $\Gamma$  with $f(\infty)=\infty$. Set $h_1:\Bbb R\to\Bbb R$ by $f\circ h_1=z$ as before. Then we have

\proclaim {Theorem 4.1} Under the above notations, the following statements are equivalent$:$

\noindent $(1)$ $\Gamma$ is a Weil-Petersson curve$;$

\noindent $(2)$ $h_1\in\WP(\Bbb R);$

\noindent $(3)$ $\arg z'\circ h^{-1}_1\in H^{\frac 12}_{\Bbb R};$

\noindent $(4)$ $h_1$ is quasisymmetric and $\arg z'\in H^{\frac 12}_{\Bbb R}$.
 \endproclaim

 \demo{Proof} From $f\circ h_1=z$ we obtain $f'=(z'\circ h^{-1}_1)(h^{-1}_1)'$, which implies that $$\Re\log f'=\log (h^{-1}_1)',\,\Im \log f'=\arg z'\circ h^{-1}_1.\tag 4.1$$
 Now $\Gamma$ is a Weil-Petersson curve if and only if
 $$\log f'\in\Cal D(\Bbb U)\Leftrightarrow \Re\log f'\in H^{\frac 12}_{\Bbb R}\Leftrightarrow\Im\log f'\in H^{\frac 12}_{\Bbb R}.$$
 By (4.1) and Proposition 2.1 we obtain that $(1)\Leftrightarrow(2)\Leftrightarrow(3)$. Now $(4)\Rightarrow(3)$ follows directly from Proposition 4.1 below. Conversely, suppose (3) holds so that (2) also holds, which implies that $h_1$ is quasisymmetric and (4) holds by Proposition 4.1 again.\quad$\square$

\enddemo

\proclaim{Proposition 4.1 ([BA], [NS]) } Let $h$ be a sense-preserving      homeomorphism  from  $\Bbb R$  onto itself. Then the pull-back operator $P_h$ defined by $P_h(u)=u\circ h$ is a bounded operator from $H^{\frac 12}$ into itself if and only if $h$ is quasisymmetric. \endproclaim

\noindent
{\bf Proof of Theorem 2.2}\quad Let $\Gamma$ be a locally rectifiable Jordan curve passing through $\infty$ and  $z=z(s)$ be an  arc-length parametrization  of  $\Gamma$. If $\Gamma$ is a Weil-Petersson curve, then it is a chord-arc curve. We conclude by David's result (see [Da]) that there exists a real-valued BMO function  $b\in \BMO_{\Bbb R}/\Bbb R$ such that $z'(s)=e^{ib(s)}$. Now Theorem 4.1 implies that $b=\arg z'\in H^{\frac 12}_{\Bbb R}/\Bbb R$. Conversely, suppose  $\Gamma$ is  a  chord-arc curve  and there exists some function $b$ in the real Sobolev class $H^{\frac 12}_{\Bbb R}/\Bbb R$  such that $z'(s)=e^{ib(s)}$. Then, as stated in section 1, $h'_1$ belongs to the class of weights $A^{\infty}$ introduced
by Muckenhoupt (see [CF], [Gar]). Thus $h_1$ is strongly quasisymmetric in the sense of Semmes [Se3] (see also [AZ], [SWe] and section 6 below) and consequently  quasisymmetric. Since $\arg z'=b\in H^{\frac 12}_{\Bbb R}/\Bbb R$, we conclude by Theorem 4.1 again that  $\Gamma$ is a Weil-Petersson curve. \quad$\square$
\vskip 0.2 cm
\noindent
{\bf Proof of Theorem 2.1}\quad Theorem 2.2 says that $\hat \Cal T=\Cal L\cap H^{\frac 12}_{\Bbb R}/\Bbb R$. Since $\Cal L$ is open in  $\BMO_{\Bbb R}/\Bbb R$ (see [Da]), we conclude by a standard discussion that
 $\hat \Cal T$ is open in  $H^{\frac 12}_{\Bbb R}/\Bbb R$ by   the continuity of the inclusion $H^{\frac 12}_{\Bbb R}/\Bbb R\hookrightarrow \BMO_{\Bbb R}/\Bbb R$.\quad$\square$
 \head 5 On Weil-Petersson quasiconformal mappings
 \endhead
 In this section, we give some preliminary results on Weil-Petersson quasiconformal mappings, which will be frequently used in the rest of the paper. They also generalize Proposition 2.1 from the real line case to the setting of Weil-Petersson curves and have independent interests of their own.

 \proclaim {Proposition 5.1} Let $F$ be a Weil-Petersson class quasiconformal mapping on the whole plane with $F(\infty)=\infty$. Then $\Gamma=F(\hat\Bbb R)$ is a Weil-Petersson curve. Furthermore, $h=F|_{\Bbb R}$ is locally absolutely continuous such that $\log h'\in H^{\frac 12}$.
 \endproclaim

\demo{Proof} Let $g$ be a Riemann mapping which takes the lower half plane $\Bbb U^*$ to the right domain $\Omega^*$ bounded by $\Gamma$. Noting that $g^{-1}\circ F$ is a Weil-Petersson  quasiconformal mapping of the lower half plane $\Bbb U^*$ onto itself, we conclude by Proposition 2.1 that $\tilde h=(g^{-1}\circ F)|_{\Bbb R}$ belongs to the Weil-Petersson class $\WP(\Bbb R)$, which implies by Proposition 2.1 again that the inverse mapping ${\tilde h}^{-1}\in \WP(\Bbb R)$ can be extended to a Weil-Petersson  quasiconformal mapping $H$ of the upper half plane $\Bbb U$ onto itself which is bi-Lipschitz under the Poincar\'e metric $\lambda_{\Bbb U}(z)|dz|$. Set $\tilde g=F\circ H$. Then $\tilde g$ is a quasiconformal extension of $g$ to the upper half plane. Now $\mu(\tilde g)=\frac{\overline{\partial}{\tilde g}}{\partial{\tilde g}}$, $\mu(F)=\frac{\overline{\partial}{F}}{\partial{F}}$ and $\mu(H)=\frac{\overline{\partial}{H}}{\partial{H}}$ satisfy
$$\mu(\tilde g)=\frac{\mu(H)+(\mu(F)\circ H)\frac{\overline{\partial H}}{\partial H}}{1-\overline{\mu(H)}(\mu(F)\circ H)\frac{\overline{\partial H}}{\partial H}},$$
which implies that
$$|\mu(\tilde g)|^2\lesssim |\mu(H)|^2+|\mu(F)\circ H|^2.$$
Since $\mu(H)\in \Cal L^{\infty}(\Bbb U)$, $\mu(F)\in \Cal L^{\infty}(\Bbb U)$, and $H$  is bi-Lipschitz under the Poincar\'e metric $\lambda_{\Bbb U}(z)|dz|$, it is easy to see that $\mu(\tilde g)\in \Cal L^{\infty}(\Bbb U)$, that is, $\tilde g$ is a Weil-Petersson quasiconformal  extension of $g$ to the whole plane, which implies by definition that $\Gamma$ is a Weil-Petersson curve. Similarly,  a Riemann mapping $f$ which takes the upper half plane $\Bbb U$ to the left domain $\Omega$ bounded by $\Gamma$ can also be extended to a Weil-Petersson quasiconformal mapping on the whole plane.

Now since $\Gamma$ is a Weil-Petersson curve, we conclude that $\log g'\in\Cal D(\Bbb U^*)$.   From $h=F|_{\Bbb R}=g\circ\tilde h$  we obtain that $h$ is locally absolutely continuous, and
$$\log h'=\log g'\circ\tilde h+\log \tilde h'=P_{\tilde h}(\log g')+\log \tilde h',$$
which implies by Propositions 2.1 and  4.1 that $\log h'\in H^{\frac 12}$ as required. \quad$\square$

\enddemo
The next result gives the converse to Proposition 5.1.

\proclaim {Proposition 5.2} Let $h$ be a sense-preserving homeomorphism from the real line onto a Weil-Petersson curve $\Gamma$ such that $h$ is locally absolutely continuous with $\log |h'|\in H^{\frac 12}$. Then $h$ can be extended to  a Weil-Petersson  quasiconformal mapping on the whole plane.
 \endproclaim
\demo{Proof} Let $z=z(s)$ be an arc-length parametrization of $\Gamma$ so that $|z'|=1$, and $f$ be a Riemann mapping from the upper half plane $\Bbb U$ onto the left domain $\Omega$ bounded by $\Gamma$ so that $\log f'\in\Cal D(\Bbb U)$. Consider two increasing homeomorphisms $h_1$ and $h_2$ of the real line $\Bbb R$ onto itself by $z\circ h_1=f$ and  $z\circ h_2=h$, respectively. Noting that $|z'|=1$, we obtain $h'_1=|f'|$, $h'_2=|h'|$, which implies that $\log h'_1=\Re\log f'\in H^{\frac 12}_{\Bbb R}$, $\log h'_2=\log |h'|\in H^{\frac 12}_{\Bbb R}$. We conclude by Proposition 2.1 that both $h_1$ and $h_2$ are in the Weil-Petersson class $\WP(\Bbb R)$. Consequently, $h^{-1}_2\circ h_1$ also belongs to the Weil-Petersson class $\WP(\Bbb R)$ and can be extended to a  Weil-Petersson  quasiconformal mapping $H$ to the upper half plane $\Bbb U$ onto itself which is bi-Lipschitz under the Poincar\'e metric $\lambda_{\Bbb U}(z)|dz|$. Then $F=f\circ H^{-1}$ is a Weil-Petersson  quasiconformal extension of $f\circ  h^{-1}_1\circ h_2=h$ from  the upper half plane $\Bbb U$ onto $\Omega$.

By the same way, we can extend $h$ to a Weil-Petersson  quasiconformal mapping  from  the lower half plane $\Bbb U^*$ onto the right domain $\Omega^*$ bounded by $\Gamma$.\quad$\square$

\enddemo

In the next section, we need to extend the arc-length parametrization of a Weil-Petersson curve to a  Weil-Petersson  quasiconformal mapping on the whole plane which is bi-Lipschitz under the Euclidian metric. Actually, we have the following general result.
\proclaim {Proposition 5.3} Let $h$ be a sense-preserving homeomorphism from the real line onto a Weil-Petersson curve $\Gamma$ such that $h$ is bi-Lipschitz under the Euclidian metric and $\log |h'|\in H^{\frac 12}$. Then $h$ can be extended to  a Weil-Petersson  quasiconformal mapping on the whole plane which is bi-Lipschitz under the Euclidian metric.
 \endproclaim

\demo{Proof}  By  Proposition 5.2 we conclude that $h$ can be extended to  a Weil-Petersson  quasiconformal mapping $F$ on the whole plane. It needs to show that $F$ is bi-Lipschitz under the Euclidian metric under the additional assumption that $h$ is bi-Lipschitz under the Euclidian metric. We only consider the  upper half-plane case. The lower half-plane case can be treated similarly. We clarify  the proof  from Semmes [Se3].

To prove the bi-Lipschitzness of $F$ under the  Euclidian metric, it is enough to show that $|\partial F|$ is bounded above and below from zero. Under the notations during the proof of Proposition 5.2,  $H|_{\Bbb R}=h^{-1}\circ z\circ h_1$, which implies that $(H|_{\Bbb R})'\asymp h'_1$ by the  bi-Lipschitzness  assumption of $h$ under the  Euclidian metric. Noting that
 $$|\partial F|=|(f'\circ H^{-1})\partial H^{-1}|=\frac{|f'||\partial H|}{|\partial H|^2-|\bar\partial H|^2}\circ H^{-1},$$
 it is sufficient to show that $|f'|$ is comparable with $|\partial H|$ on the upper half plane.

 Since $H$ is bi-Lipschitz under the Poincar\'e metric $\lambda_{\Bbb U}(z)|dz|$, for $z=x+iy\in\Bbb H$ we conclude by the distortion theorem for quasiconformal mappings that
 $$
 |\partial H(z)|\asymp \frac{\Im H(z)}{\Im z}\asymp\frac{H(x+y)-H(x-y)}{y}=\frac{1}{y}\int_{x-y}^{x+y}H'(t)dt\asymp\frac{1}{y}\int_{x-y}^{x+y}h'_1(t)dt. \tag 5.1$$
 Since $h_1$ belongs to the Weil-Petersson class $\WP(\Bbb R)$, which implies that $h'_1$ belongs to the class of weights $A^{\infty}$ introduced
by Muckenhoupt (see [CF], [Gar]), we obtain that
$$
\frac{1}{2y}\int_{x-y}^{x+y}h'_1(t)dt\asymp\exp\left(\frac{1}{2y}\int_{x-y}^{x+y}\log h'_1(t)dt\right).\tag 5.2$$

To estimate the right side of (5.2), we need a well-known result on BMO functions (see [Gar]): For $u\in\BMO$, it holds that
$$\sup_{z\in\Bbb U}\int_{\Bbb R}|u(t)-u(z)|P_z(t)dt\asymp\|u\|_{\BMO},\tag 5.3$$
where
$$P_z(t)=\frac{1}{\pi}\frac{y}{(x-t)^2+y^2}$$
is the Poisson kernel, while $u(z)$ is the Poisson integral of $u(x)$, namely,
$$u(z)=\frac{1}{\pi}\int_{\Bbb R}\frac{y}{(x-t)^2+y^2}u(t)dt.$$
Thus, for $z=x+iy\in\Bbb U$ we have
$$\align
\left|\frac{1}{2y}\int_{x-y}^{x+y}u(t)dt-u(z)\right|&=\left|\frac{1}{2y}\int_{x-y}^{x+y}(u(t)-u(z))dt\right|\\
&\le\frac{1}{2y}\int_{x-y}^{x+y}|u(t)-u(z)|dt\tag 5.4\\
&\le\int_{x-y}^{x+y}|u(t)-u(z)|P_z(t)dt\lesssim\|u\|_{\BMO}.
\endalign$$
Since $\log f'\in\Cal D(\Bbb U)$, $\log |f'(z)|$ is the Poisson integral of $\log |f'(x)|=\log h'_1(x)$, which implies by (5.4) that
$$\left|\frac{1}{2y}\int_{x-y}^{x+y}\log h'_1(t)dt-\log |f'(z)|\right|\lesssim\|\log h'_1\|_{\BMO}\lesssim\|\log h'_1\|_{H^{\frac 12}}.\tag 5.5$$
We conclude from (5.1), (5.2) and (5.5) that $|f'|$ is comparable with $|\partial H|$ as desired.\quad$\square$

\enddemo

\proclaim {Corollary 5.4} Let $z=z(s)$ be an arc-length parametrization of a Weil-Petersson curve $\Gamma$. Then $z$ can be extended to  a Weil-Petersson  quasiconformal mapping on the whole plane which is bi-Lipschitz under the Euclidian metric.
 \endproclaim
\demo{Proof} In Proposition 5.3, replacing  $h$ with the arc-length parametrization $z$, we conclude that $z$ can be extended to  a Weil-Petersson  quasiconformal mapping on the whole plane which is bi-Lipschitz under the Euclidian metric. \quad$\square$

\enddemo
 \head 6 More on the openness of $\hat\Cal T$
\endhead
In section 4, we  proved that $\hat\Cal T$ is an open subset of  $H^{\frac 12}_{\Bbb R}/\Bbb R$ by means of the openness of $\Cal L$, which  implies that the set of all normalized Weil-Petersson curves, a model of the Weil-Petersson Teichm\"uller space $\Cal T$, can be  endowed with  a real Hilbert manifold structure in a geometric manner by the correspondence $\Gamma\mapsto b$. However, to prove Theorems 2.3 and 2.4, which says that  this new real Hilbert manifold structure is topologically equivalent to the  standard complex Hilbert manifold structure  given by Takhtajan-Teo [TT], we need more information about the the openness of $\hat\Cal T$ in $H^{\frac 12}_{\Bbb R}/\Bbb R$. For any $b\in\hat\Cal T$, we will show that there exists some neighbourhood $U(b)$ in  $H^{\frac 12}/\Bbb C$ such that for each $u\in U(b)$ the induced mapping $\gamma_u$ defined by (2.11)  can be extended to a  Weil-Petersson  quasiconformal on the whole plane whose Beltrami coefficient depends holomorphically on $u$ (see Proposition 6.2 and Theorem 6.1 below), a fact which will play an essential role in the proof of Theorems 2.3 and 2.4. In particular, $\gamma_u$ maps the real line $\hat\Bbb R$ onto a Weil-Ptersson curve and $b\in\hat\Cal T$ is an interior point.

We first quickly review some results  in our paper [ST], where we explored such an approach at the base point $0\in\hat\Cal T$.  We begin  with a basic result of Coifman-Meyer [CM]. For $u\in\BMO$ on the real line,  set as before that
$$\gamma_u(x)=\int_0^xe^{iu(t)}dt,\quad x\in\Bbb R.\tag 6.1$$
Coifman-Meyer [CM] showed that $\gamma_u$ is a strongly quasisymmetric homeomorphism from the extended real line $\hat\Bbb R$ onto a chord-arc curve $\Gamma_u=\gamma_u(\hat\Bbb R)$ when $\|u\|_{\BMO}$ is small. Here   a sense preserving homeomorphism $h$ on $\hat\Bbb R$ is said to be
strongly quasisymmetric if it is locally absolutely continuous so that $|h'|\in A^{\infty}$ and it maps $\hat\Bbb R$ onto a chord-arc curve passing through the point at infinity (see [Se3]). Later, Semmes [Se3] showed that, when $\|u\|_{\BMO}$ is small,  $\gamma_u$ can be extended a quasiconformal mapping to the whole plane whose Beltrami coefficient satisfies certain Carleson measure condition. To be precise, let $\varphi$ and $\psi$ be two $C^{\infty}$ real-valued function on the real line  supported on $[-1, 1]$ such that $\varphi$ is even, $\psi$ is odd and $\int_{\Bbb R}\varphi(x)dx=1$, $\int_{\Bbb R}\psi(x)xdx=1$. Define
$$\rho(x, y)=\rho_u(x, y)=\varphi_y\ast \gamma_u(x)-i(\text{sgn} y)\psi_y\ast \gamma_u(x),\quad  z=x+iy\in\Bbb U\cup\Bbb U^*,\tag 6.2$$
and $\rho(x, 0)=\gamma_u(x)$ for $x\in\Bbb R$.
Then $\rho$ is a quasiconformal mapping on the whole plane whose Beltrami coefficient  satisfies certain Carleson measure condition when    $\|u\|_{\BMO}$ is small.
We proved in [ST]  that $\rho_u$ is in the  Weil-Petersson class   when    $u\in H^{\frac 12}$ is small.

\proclaim{Proposition  6.1 ([ST])} There exists some universal constant $\delta>0$ such that,  for any  $u\in U(0, \delta)\doteq\{u\in H^{\frac 12}/\Bbb C: \|u\|_{H^{\frac 12}}<\delta\}$, the mapping $\rho=\rho_u$ defined by (6.2) is a Weil-Petersson  quasiconformal extension of $\gamma_u$ on the whole plane whose Beltrami coefficient $\mu$ satisfies $\|\mu|_{\Bbb U}\|_{\WP}\lesssim\|u\|_{H^{\frac 12}}$ and $\|\mu|_{\Bbb U^*}\|_{\WP}\lesssim\|u\|_{H^{\frac 12}}$.
\endproclaim

By Proposition 5.1, we conclude that $\Gamma_u=\gamma_u(\hat\Bbb R)$ is a Weil-Petersson curve when $\|u\|_{H^{\frac 12}}$ is small. Moreover, when $u\in H_{\Bbb R}^{\frac 12}/\Bbb R$, $\gamma_u$ is the normalized arc-parametrization $z_u$ of the normalized Weil-Petersson curve $\Gamma_u$. Consequently, when $u\in H_{\Bbb R}^{\frac 12}/\Bbb R$ is small, $u\in\hat\Cal T$, which implies that $0$ is an interior point of $\hat\Cal T$.

During the proof of Proposition 6.1, we established the following result, which will be frequently used later.
\proclaim{Lemma 6.1 ([ST])} Let $\mu\in L^{\infty}(\Bbb U)$ and $u\in H^{\frac 12}$ satisfy the following condition
$$|\mu(x+iy)|^2\lesssim \frac{1}{y}\int^y_{-y}|u(t+x)-u(x)|^2dt.\tag 6.3$$
Then $\mu\in \Cal L^{\infty}(\Bbb U)$.
\endproclaim


We also recall the following result from [ST], which will be used in the proof of Theorem 2.3 in the next section.

\proclaim{Proposition 6.2 ([ST])}  For $u\in U(0, \delta)$, let $\Lambda(u)$ denote the Beltrami coefficient on the upper half plane $\Bbb U$ for the quasiconformal mapping $\rho_u$ defined by (6.2). Then $\Lambda: U(0, \delta)\to\Cal M(\Bbb U)$ is holomorphic.
\endproclaim

In the rest of this section, we will extend the above approach to a general point of $\hat\Cal T$. Let $b\in \hat\Cal T$ be a non-zero element and $\Gamma_b$ be the  normalized Weil-Petersson curve whose normalized arc-length parametrization $z_b$ satisfies $z'_b=e^{ib}$. By Corollary 5.4, there exists   a Weil-Petersson  quasiconformal mapping $\tau$ on the whole plane which is bi-Lipschitz under the Euclidian metric and satisfies $\tau(x)=z_b(x)$ for $x\in\Bbb R$. Now for any $u\in H^{\frac 12}$,  set
$$\omega_u(x)=\int_0^xe^{i(b(t)+u(t))}dt,\quad x\in\Bbb R.\tag 6.5$$
We will show that $\Gamma_{b+u}=\gamma_{b+u}(\hat\Bbb R)=\omega_u(\hat\Bbb R)$ is a Weil-Petersson curve when $u\in H^{\frac 12}$ is small. In particular, this  implies that $b+u\in\hat\Cal T$ when $u\in H^{\frac 12}_{\Bbb R}/\Bbb R$ is small, that is, $b\in \hat\Cal T$ is an interior point. By Proposition 5.1, it is sufficient  to show that $\omega_u$ can be extended to a Weil-Petersson quasiconformal mapping on the whole plane when $u\in H^{\frac 12}$ is small.

We use Semmes' construction (see [Se3]). Let $\varphi$  be a $C^{\infty}$ real-valued even function on the real line  supported on $[-1, 1]$ such that $\int_{\Bbb R}\varphi(x)dx=1$.
On setting  $\varphi_y(x)=|y|^{-1}\varphi(|y|^{-1}x)$ for $y\neq 0$, we also want that  $|\varphi_y\ast z'_b|\ge \epsilon_0$ for some  $\epsilon_0>0$. As pointed out by Semmes [Se3], this can be done as soon as $\Gamma_b$ is a chord-arc curve, especially a Weil-Petersson curve. As in Lemma 3.1, we consider
 $$R_y(w)(x)=\frac{\varphi_y\ast (z'_bw)(x)}{\varphi_y\ast z'_b(x)}. \tag 6.6$$
Now we are ready to define
$$\rho(z)=\rho_u(z)=\varphi_{y}\ast \omega_u(x)+R_y(e^u)(x)\{\tau(z)-\varphi_y\ast\tau(x)\},\,z=x+iy\in\Bbb U\cup\Bbb U^*,\tag 6.7$$
and $\rho(x)=\omega_u(x)$ for $x\in\Bbb R$.
Semmes [Se3] showed that $\rho$ is quasiconformal on the whole plane
when $\|u\|_{\BMO}$ is small. We will show that $\rho$ is in the Weil-Petersson class on the whole plane when $\|u\|_{H^{\frac 12}}$ is small.
We only consider the  upper half-plane case. The lower half-plane case can be treated similarly.

We proceed to estimate the derivatives of $\rho$. From (6.7) we have
$$
\align
\bar{\partial}\rho(z)
=&R_y(e^u)(x)\bar\partial\tau(z)+\bar\partial (R_y(e^u)(x))\{\tau(z)-\varphi_y\ast\tau(x)\}\\
&+\bar{\partial}(\varphi_{y}\ast \omega_u(x))-R_y(e^u)(x)\bar{\partial}(\varphi_y\ast\tau(x)). \tag 6.8
\endalign
$$
For $x\in\Bbb R$ and $y>0$, consider $I=[x-y, x+y]$ as before so that
 $$u_I=\frac{1}{2y}\int^{x+y}_{x-y}u(t)dt.$$
Since $\int_{\Bbb R}\varphi(x)dx=1$, which implies that $\int_{\Bbb R}\varphi_y(x)dx=1$,  we obtain
$$|\tau(z)-\varphi_y\ast\tau(x)|=|\varphi_y\ast(\tau-\tau(z))(x)|\lesssim\frac{1}{|I|}\int_I|\tau(t)-\tau(z)|dt\lesssim|z-t|\lesssim y,\tag 6.9$$
since $\tau$ is bi-Lipschitz under the Euclidean metric.

Now set
$$\psi(x)=\frac 12((1-ix)\varphi(x))'=\frac 12(\varphi'(x)-i(\varphi(x)+x\varphi'(x))).$$
Clearly, $\psi$ is a  $C^{\infty}$ function  on the real line which is supported on $[-1, 1]$ and satisfy $\int_{\Bbb R}\psi(x)dx=0$.
A direct computation yields that
$$y\bar\partial{(\varphi_{y}\ast w(x))}=\psi_y\ast w(x). $$
Since $|\varphi_y\ast z'_b|\ge \epsilon_0$, we have
$$
\align
&|y\bar\partial (R_y(e^u)(x))|\\
&=\left|\frac{y\bar\partial(\varphi_{y}\ast ({z_b}'e^u)(x))-yR_y(e^u)(x)\bar\partial(\varphi_{y}\ast z'_b(x))}{\varphi_{y}\ast z'_b(x)}\right|\\
&=\left|\frac{\psi_y\ast (z'_be^u)(x)-R_y(e^u)(x)\psi_y\ast z'_b(x)}{\varphi_{y}\ast z'_b(x)}\right|\\
&\lesssim|\psi_y\ast (z'_b(e^u-e^{u_I}))(x)|+|\psi_y\ast z'_b(x)||e^{u_I}-R_y(e^u)(x)|. \tag 6.10
\endalign
$$

For the first term in (6.10),
$$
\align
|\psi_y\ast (z'_b(e^u-e^{u_I}))(x)|&\lesssim \frac{1}{|I|}\int_I|z'_b(t)(e^{u(t)}-e^{u_I})|dt\\
&\lesssim \frac{1}{|I|}\int_I|e^{u(t)-u_I}-1||e^{u_I}|dt\\
&\lesssim \frac{1}{|I|}\int_I|u(t)-u_I||e^{u(t)}|dt.\tag 6.11
\endalign
$$
Then, using H\"older inequality we obtain from  Lemma 3.1 and  (3.4), (3.5), (3.8) that
$$
\align
\frac{|\psi_y\ast (z'_b(e^u-e^{u_I}))(x)|}{|R_y(e^u)(x)|}&\lesssim \frac{1}{|I|}\int_I|u(t)-u_I||e^{u(t)-R_y(u)(x)}|dt\\
&\lesssim \frac{1}{|I|}\int_I|u(t)-u_I|e^{|u(t)-u_I|+|R_y(u)(x)-u_I|}dt\\
&\lesssim\|u\|_{\BMO}\tag 6.12
\endalign
$$
when $\|u\|_{\BMO}$ is small.

For the second term in (6.10),
$$
|\psi_y\ast z'_b(x)||e^{u_I}-R_y(e^u)(x)|\lesssim |R_y(e^u-e^{u_I})(x)|
\lesssim |\varphi_y\ast (z'_b(e^u-e^{u_I}))(x)|.
$$
Repeat the reasoning in (6.11) and (6.12), we have
$$\frac{|\psi_y\ast z'_b(x)||e^{u_I}-R_y(e^u)(x)|}{|R_y(e^u)(x)|}\lesssim \|u\|_{\BMO}\tag 6.13
$$
when $\|u\|_{\BMO}$ is small.
So by (6.9), (6.10), (6.12) and (6.13) we have
$$|\bar\partial (R_y(e^u)(x))\{\tau(z)-\varphi_y\ast \tau(x)\}|\lesssim  y|\bar\partial (R_y(e^u)(x))|\lesssim \|u\|_{\BMO} |R_y(e^u)(x)| \tag 6.14$$
when $\|u\|_{\BMO}$ is small.

Next, we will prove $|y\bar\partial (R_y(e^u)(x))||R_y(e^u)(x)|^{-1}$ is square integrable in the Poincar\'{e} metric when $\|u\|_{ H^{\frac 12}}$ and consequently $\|u\|_{\BMO}$ is small.
Similar to (6.10),
$$
\align
&|y\bar\partial (R_y(e^u)(x))|\\
&=\left|\frac{\psi_y\ast (z'_be^u)(x)-R_y(e^u)(x)\psi_y\ast z'_b(x)}{\varphi_{y}\ast z'_b(x)}\right|\\
&\lesssim|\psi_y\ast (z'_b(e^u-e^{u(x)}))(x)|+|\psi_y\ast z'_b(x)||e^{u(x)}-R_y(e^u)(x)|. \tag 6.15
\endalign
$$

For the first part of (6.15), similar to (6.11)
$$
|\psi_y\ast (z'_b(e^u-e^{u(x)}))(x)|\lesssim\frac{1}{|I|}\int_I|u(t)-u(x)||e^{u(t)}|dt.\tag 6.16
$$
Then, by lemma 3.1,
$$
\frac{|\psi_y\ast (z'_b(e^u-e^{u(x)}))(x)|}{|R_y(e^u)(x)|}\lesssim \frac{1}{|I|}\int_I|u(t)-u(x)||e^{u(t)-R_y(u)(x)}|dt.\tag 6.17
$$
By H\"{o}lder inequality and (3.9),  we conclude
$$
\align
\frac{|\psi_y\ast (z'_b(e^u-e^{u(x)}))(x)|^2}{|R_y(e^u)(x)|^2}&\lesssim \frac{1}{|I|^2}\int_I|u(t)-u(x)|^2dt\int_I|e^{u(t)-R_y(u)(x)}|^2dt\\
&\lesssim\frac{1}{|I|}\int_I|u(t)-u(x)|^2dt\\
&\lesssim\frac{1}{y}\int^{y}_{-y}|u(t+x)-u(x)|^2dt.\tag 6.18
\endalign$$
Consequently, by Lemma 6.1 and (6.4),
$$\iint_{\Bbb U}\frac{|\psi_y\ast (z'_b(e^u-e^{u(x)}))(x)|^2}{|R_y(e^u)(x)|^2}\frac{1}{|y|^2}dxdy\lesssim ||u||_{H^{\frac{1}{2}}}^2.\tag 6.19$$

For the second part of (6.15), $$|\psi_y\ast z'_b||e^{u(x)}-R_y(e^u)(x)|\lesssim |R_{y}(e^{u}-e^{u(x)})(x)|\lesssim |\varphi_{y}\ast (z'_b(e^{u}-e^{u(x)}))(x)|.$$
Doing the same as (6.16)-(6.19), we can obtain that
$$\iint_{\Bbb U}\frac{|\varphi_{y}\ast (z'_b(e^{u}-e^{u(x)}))(x)|^2}{|R_y(e^u)(x)|^2}\frac{1}{|y|^2}dxdy\lesssim \|u\|_{H^{\frac{1}{2}}}^2.\tag 6.20 $$
Therefore, by (6.15), (6.19) and (6.20) we obtain that
$$\iint_{\Bbb U}\frac{|\bar\partial (R_y(e^u)(x))|^2}{|R_y(e^u)(x)|^2}dxdy\lesssim \|u\|_{H^{\frac{1}{2}}}^2\tag 6.21$$
when $\|u\|_{H^{\frac 12}}$ is small.

Now, we consider the third and forth parts of (6.8). We have
$$
\align
&\frac{\partial}{\partial x}(\varphi_{y}\ast \omega_u(x))-R_y(e^u)(x)\frac{\partial}{\partial x}(\varphi_y\ast\tau(x))\\
&=
\varphi_{y}\ast \omega'_u(x)-R_y(e^u)(x)\varphi_y\ast\tau'(x)\tag 6.22\\
&=\varphi_{y}\ast (z'_be^{iu})(x)-R_y(e^u)(x)\varphi_y\ast z'_b(x)=0.
\endalign$$
Noting that
$$\frac{\partial}{\partial y}(\varphi_{y}\ast w(x))=\alpha_y\ast w'(x),$$
where $\alpha(x)=-x\varphi(x)$ is a  $C^{\infty}$ function on the real line  which is supported on [-1,1] and satisfies $\int_{\Bbb R}\alpha(x)dx=0$, we have
$$\frac{\partial}{\partial y}(\varphi_{y}\ast \omega_u(x))-R_y(e^u)(x)\frac{\partial}{\partial y}(\varphi_y\ast\tau(x))=\alpha_y\ast (z'_be^{u})(x)-R_y(e^u)(x)\alpha_y\ast z'_b(x).\tag 6.23$$
Noting that
$$
|\alpha_y\ast (z'_be^{u})(x)-R_y(e^u)(x)\alpha_y\ast z'_b(x)|\le|\alpha_y\ast (z'_b(e^{u}-e^{u_I}))(x)|+|\alpha_y\ast (z'_b(R_y(e^u)(x)-e^{u_I}))(x)|,
$$
we do the same as (6.11), (6.12) and obtain
$$|\alpha_y\ast (z'_be^{u})(x)-R_y(e^u)(x)\alpha_y\ast z'_b(x)|\lesssim|R_{y}(e^u)(x)|\|u\|_{\BMO}\tag 6.24$$
when $\|u\|_{\BMO}$ is small.
Similar to (6.16)-(6.19), we can prove that
$$\iint_{\Bbb U}\frac{|\alpha_y\ast (z'_be^{u})(x)-R_y(e^u)(x)\alpha_y\ast z'_b(x)|^2}{|R_y(e^u)(x)|^2}\frac{1}{y^2}dxdy\lesssim \|u\|_{H^{\frac{1}{2}}}^2\tag 6.25$$
when $\|u\|_{H^{\frac 12}}$ is small.

Summarizing the above,  we have
$$|\bar\partial\rho(z)-R_y(e^u)(x)\bar\partial{\tau(z)}|\lesssim|R_{y}(e^u)(x)|||u||_{H^{\frac{1}{2}}},\tag 6.26 $$
and
$$\iint_{\Bbb U}\frac{|\bar\partial\rho(z)-R_y(e^u)(x)\bar\partial{\tau(z)}|^2}{|R_{y}(e^u)(x)|^2}\frac{1}{y^2}dxdy\lesssim\|u\|_{H^{\frac{1}{2}}}\tag 6.27 $$
if $\|u\|_{H^{\frac{1}{2}}}$ is small enough.

For another derivative of $\rho$,
$$
\align
{\partial}\rho(z)
=&R_y(e^u)(x)\partial\tau(z)+\partial (R_y(e^u)(x))\{\tau(z)-\varphi_y\ast\tau(x)\}\\
&+{\partial}(\varphi_{y}\ast \omega_u(x))-R_y(e^u)(x){\partial}(\varphi_y\ast\tau(x)). \tag 6.28
\endalign
$$
Similarly, we can prove that
$$|\partial\rho(z)-R_y(e^u)(x)\partial{\tau(z)}|\lesssim|R_{y}(e^u)(x)|||u||_{H^{\frac{1}{2}}},\tag 6.29 $$
and
$$\iint_{\Bbb U}\frac{|\partial\rho(z)-R_y(e^u)(x)\partial{\tau(z)}|^2}{|R_{y}(e^u)(x)|^2}\frac{1}{y^2}dxdy\lesssim\|u\|_{H^{\frac{1}{2}}}\tag 6.30 $$
if $\|u\|_{H^{\frac{1}{2}}}$ is small enough.
Since  $\tau$ is bi-Lipschitz under the Euclidean metric, which implies that $|\partial\tau|\asymp 1$, we obtain from (6.26-6.30) that
$$\left\|\frac{\bar\partial\rho}{\partial\rho}-\frac{\bar\partial\tau}{\partial\tau}\right\|_{\infty}\lesssim\|u\|_{H^{\frac 12}},\tag 6.31$$
and
$$\iint_{\Bbb U}\left|\frac{\bar\partial\rho}{\partial\rho}(z)-\frac{\bar\partial\tau}{\partial\tau}(z)\right|^2\frac{1}{y^2}dxdy\lesssim\|u\|_{H^{\frac 12}}\tag 6.32$$
when $\|u\|_{H^{\frac 12}}$ is small.
In particular,  $\rho$ is in the Weil-Petersson class on the upper half plane when $\|u\|_{H^{\frac 12}}$ is small. \quad$\square$

We summarize the above discussion in the following

\proclaim{Theorem 6.1} For each non-zero $b\in\hat\Cal T$, there exists some  $\delta>0$  so that for each $v\in U(b, \delta)\doteq \{v\in H^{\frac 12}/\Bbb C: \|v-b\|_{H^{\frac 12}}<\delta\}$, the mapping $\rho_{v-b}$ defined by (6.7) is a Weil-Petersson quasiconformal extension of $\gamma_v$  on the whole plane.  For $v\in U(b, \delta)$, let $\Lambda(v)$ denote the Beltrami coefficient for the quasiconformal mapping $\rho_{v-b}$ on the upper half plane. Then $\Lambda: U(b, \delta)\to\Cal M(\Bbb U)$ is holomorphic.
\endproclaim

 \demo{Proof} It remains to show the holomorphy of $\Lambda$. The proof is almost the same as the one of Proposition 6.2 given in [ST] and is reproduced here for the sake of completeness. By (6.31) and (6.32), $\Lambda$ is bounded in $U(b, \delta)$. So it is sufficient to show that, for each fixed pair of $(u, v)$ with $u\in U(b, \delta)$, $v\in H^{\frac 12}/\Bbb C$, $\tilde\Lambda(t)\doteq\Lambda(u+tv)$ is holomorphic in a small neighbourhood of $t=0$ in the complex plane. To do so, choose
 $$0<\epsilon<\frac{\delta-\|u-b\|_{H^{\frac 12}}}{2\|v\|_{H^{\frac 12}}}$$
so that $u+tv\in U(b, \delta)$ when $|t|\le 2\epsilon$. We conclude by (6.8) and (6.28) that  $\tilde\Lambda(t)(z)$ is holomorphic in $|t|\le 2\epsilon$ for fixed $z\in\Bbb U$. For $|t_0|<\epsilon$, $|t|<\epsilon$, Cauchy formula yields that
$$\align
\left|\frac{\tilde\Lambda(t)(z)-\tilde\Lambda(t_0)(z)}{t-t_0}-\frac{d}{dt}|_{t=t_0}\tilde\Lambda(t)(z)\right|
&=\frac{|t-t_0|}{2\pi }\left|\int_{|\zeta|=2\epsilon}\frac{\tilde\Lambda(\zeta)(z)}{(\zeta-t)(\zeta-t_0)^2}d\zeta\right|\\
&\le\frac{|t-t_0|}{2\pi\epsilon^3}\int_{|\zeta|=2\epsilon}|\tilde\Lambda(\zeta)(z)||d\zeta|.
\endalign
$$
Thus, by (6.31),
$$
\left\|\frac{\tilde\Lambda(t)-\tilde\Lambda(t_0)}{t-t_0}-\frac{d}{dt}|_{t=t_0}\tilde\Lambda(t)\right\|_{\infty}
\le\frac{|t-t_0|}{2\pi\epsilon^3}\int_{|\zeta|=2\epsilon}\|\tilde\Lambda(\zeta)\|_{\infty}|d\zeta|\le C(u, v)|t-t_0|,
$$
and by (6.32),
$$\align
&\iint_{\Bbb U}\frac{1}{y^2}\left|\frac{\tilde\Lambda(t)(z)-\tilde\Lambda(t_0)(z)}{t-t_0}-\frac{d}{dt}|_{t=t_0}\tilde\Lambda(t)(z)\right|^2dxdy\\
&\le\frac{|t-t_0|^2}{4\pi^2\epsilon^6}\iint_{\Bbb U}\frac{1}{y^2}\left(\int_{|\zeta|=2\epsilon}|\tilde\Lambda(\zeta)(z)||d\zeta|\right)^2dxdy\\
&\le\frac{|t-t_0|^2}{\pi\epsilon^5}\iint_{\Bbb U}\int_{|\zeta|=2\epsilon}\frac{|\tilde\Lambda(\zeta)(z)|^2}{y^2}|d\zeta|dxdy\\
&=\frac{|t-t_0|^2}{\pi\epsilon^5}\int_{|\zeta|=2\epsilon}\iint_{\Bbb U}\frac{|\tilde\Lambda(\zeta)(z)|^2}{y^2}dxdy|d\zeta|\\
&\lesssim C(u, v)|t-t_0|^2.
\endalign
$$
 Consequently, the limit
$$\lim_{t\to t_0}\frac{\tilde\Lambda(t)-\tilde\Lambda(t_0)}{t-t_0}=\frac{d}{dt}|_{t=t_0}\tilde\Lambda(t)$$
exists in $\Cal M(\Bbb U)$ and $\Lambda: U(b, \delta)\to\Cal M(\Bbb U)$ is holomorphic. \quad$\square$
\enddemo

\vskip 0.2 cm
\noindent
{\bf Remark.}\quad Recall that  $\hat\Cal T_e$ is the set of all $u\in H^{\frac 12}/\Bbb C$ such that $\gamma_u$ is a homeomorphism from the real line $\Bbb R$ onto a Weil-Petersson curve, and  $\hat\Cal T=\hat\Cal T_e\cap H_{\Bbb R}^{\frac 12}/\Bbb R$. Our discussion not only shows that $\hat \Cal T$ is an open set in $H_{\Bbb R}^{\frac 12}/\Bbb R$, namely, each point $b\in \hat \Cal T$ is an interior point of $\hat \Cal T$, but also shows that  each point $b\in \hat \Cal T$ is actually an interior point of $\hat \Cal T_e$. In fact,  $\hat\Cal T_e$ is  an open subset of $H^{\frac 12}/\Bbb C$. We will come back this in the last section.

\head 7 Proof of Theorems 2.3 and 2.4
\endhead

For $b\in\hat\Cal T$, we assume without loss of generality that $\int_0^1e^{ib(t)}dt>0$, and set as before
$$z_b(s)=\int_0^se^{ib(t)}dt.\tag 7.1$$
We also use the notations  $\Gamma_b=z_b(\hat\Bbb R)$,  $f_b$, $g_b$, and $h_b$ introduced in section 2. Now we consider $h_1=f_b^{-1}\circ z_b$ and $h_2=g_b^{-1}\circ z_b$. Then we have the following result.

\proclaim{Theorem 7.1} Both $h_1$ and $h_2$ depend on $b$ real-analytically. Precisely, the correspondence $b\mapsto h_1$ induces a real-analytic map  from $\hat\Cal T$ into the normalized Weil-Petersson class $\WP_0(\Bbb R)\,(=\Cal T)$, and so does the correspondence $b\mapsto h_2$.
\endproclaim
\demo{Proof} For each $b\in\hat\Cal T$, we consider the neighborhood $U(b, \delta)$ in Proposition 6.1 or Theorem 6.1. For each $v\in U(b, \delta)$, we denote as above by $\Lambda(v)$ the Beltrami coefficient on the upper half plane $\Bbb U$ of the quasiconformal mapping $\rho_{v-b}$ defined by (6.2) or (6.7). Then $\Lambda: U(b, \delta)\to\Cal M(\Bbb U)$ is holomorphic by Proposition 6.2 or Theorem 6.1, which implies that $\Lambda$ is real-analytic from $U_{\Bbb R}(b,\delta)$, the real-valued functions in $U(b, \delta)$, into $\Cal M(\Bbb U)$. On the other hand, when $v\in U_{\Bbb R}(b, \delta)$, $f_v\circ h_1=z_v=\omega_{v-b}=\rho_{v-b}|_{\Bbb R}$, which implies that $h_1=(f^{-1}_v\circ \rho_{v-b})|_{\Bbb R}$, or equivalently, $h_1=f^{\Lambda(v)}|_{\Bbb R}$. Consequently, the correspondence $b\mapsto h_1$ induces a real-analytic map  from $\hat\Cal T$ into the normalized Weil-Petersson class $\WP_0(\Bbb R)\,(=\Cal T)$. By the same way, we can prove that  the correspondence $b\mapsto h_2$ also induces a real-analytic map  from $\hat\Cal T$ into the normalized Weil-Petersson class $\WP_0(\Bbb R)\,(=\Cal T)$.\quad$\square$
\enddemo
\vskip 0.2 cm
\noindent
{\bf Proof of Theorems 2.3 and 2.4}\quad For each $b\in\hat\Cal T$, we have $h_b=h_1\circ h^{-1}_2$. Since the Weil-Petersson Teichm\"uller space $\WP_0(\Bbb R)$ is a topological group, we conclude by Theorem 7.1 that $h_b$ depends continuously on $b$, or equivalently, $g_b$ depends continuously on $b$.

It is easy to see that the correspondence $b\mapsto h_b$ induces a one-to-one map from $\hat\Cal T$ onto the normalized Weil-Petersson class $\WP_0(\Bbb R)\,(=\Cal T)$, or equivalently, the correspondence $b\mapsto \log h'_b$ induces a one-to-one map from $\hat\Cal T$ onto the real Sobolev space $H_{\Bbb R}^{\frac 12}/\Bbb R$. Thus, for each $h\in\WP_0(\Bbb R)$ there exists unique $b\in\hat\Cal T$ such that $h=h_b=f^{-1}_b\circ g_b$. Suppose $h_{b_n}\to h_b$ in $\WP_0(\Bbb R)$, or equivalently, $\|\log h'_{b_n}-\log h'_b\|_{H^{\frac 12}}\to 0$. We need to show that $\|b_n-b\|_{H^{\frac 12}}\to 0$. Writing $h_{b_n}=f^{-1}_{b_n}\circ g_{b_n}$, we have $\|\log g'_{b_n}-\log g'_{b}\|_{\Cal D(\Bbb U^*)}\to 0$. Let $z_b$ and $z_{b_n}$ denote the arc-length parametrization of the normalized Weil-Petersson curve $\Gamma_{b}$ and $\Gamma_{b_n}$, respectively. Set as above that $h_2=g_b^{-1}\circ z_b$, and $h_{2n}=g_{b_n}^{-1}\circ z_{b_n}$. Then $(h^{-1}_2)'=|g'_b|$, and $(h^{-1}_{2n})'=|g'_{b_n}|$. Noting that
$$\log (h^{-1}_{2n})'-\log  (h^{-1}_2)'=\Re(\log g'_{b_n}-\log g'_b),$$
we conclude that $\|\log (h^{-1}_{2n})'-\log  (h^{-1}_2)'\|_{H^{\frac 12}}\to 0$, that is, $h^{-1}_{2n}\to h^{-1}_2$ in $\WP_0(\Bbb R)$, or equivalently, $h_{2n}\to h_2$ in $\WP_0(\Bbb R)$.  On the other hand, from
$$z_b=g_b\circ h_2\Rightarrow z'_b=(g'_b\circ h_2)h'_2\Rightarrow ib=\log z'_b=\log(g'_b\circ h_2)+\log h'_2$$
we obtain $b=\Im\log(g'_b\circ h_2)$. Similarly, $b_n=\Im\log(g'_{b_n}\circ h_{2n})$. So we have
$$b_n-b=\Im((\log g'_{b_n})\circ h_{2n}-(\log g'_b)\circ h_2).$$
Noting that $\|\log g'_{b_n}-\log g'_{b}\|_{\Cal D(\Bbb U^*)}\to 0$, and $h_{2n}\to h_2$ in $\WP_0(\Bbb R)$, we conclude by the following result (see  Corollary 4.2 in [HWS] and also also Lemma 7.2 in [Sh]) that $\|b_n-b\|_{H^{\frac 12}}\to 0$ as required.

\proclaim{Proposition 7.1 ([HWS], [Sh])} Let $h_t$, $t\in [0, t_0]$, be quasisymmetric homeomorphisms on the real line which keep the points $0$ and $1$ fixed. Suppose $u_t:[0, t_0]\to H^{\frac12}$ and $h_t:[0, t_0]\to T$ are continuous. Then $P_{h_t}u_t: [0, t_0]\to H^{\frac 12}$ is continuous.
\endproclaim

To complete the proof, we need to show that the correspondence $b\mapsto \log f'_b$ induces a homeomorphism from $\hat\Cal T$ onto its image in $\Cal D(\Bbb U)/\Bbb C$. This can be obtained by means of the following facts: Let $J(z)=\bar z$ denote the standard reflection with respect to the real line. Then for each $b\in\hat\Cal T$, we have $z_{-b}=Jz_b$, $\Gamma_{-b}=J(\Gamma_b)$, $f_{-b}=Jg_bJ$, $g_{-b}=Jf_bJ$, $h_{-b}=h^{-1}_b$.\quad$\square$

\head 8 Generalized Weil-Petersson homeomorphisms and Proof of Theorem 2.5
\endhead
A sense-preserving homeomorphism $h$ on the real line $\Bbb R$ with $h(\infty)=\infty$ is called a generalized Weil-Petersson homeomorphism if $h$ is locally absolutely continuous with  $\log h'\in H^{\frac 12}$ (or equivalently, $\log |h'|\in H^{\frac 12}$ by Propositions 5.1 and 5.2)  and $h(\hat\Bbb R)$ is a Weil-Petersson curve. By Propositions 5.1 and 5.2, a sense-preserving homeomorphism $h$ on the real line $\Bbb R$ is  a generalized Weil-Petersson homeomorphism if and only if $h$ can be extended to a Weil-Petersson quasiconformal mapping on the whole plane with $\infty$ fixed. There are several ways to parameterize the class $\WP(\Bbb C)$  of all generalized Weil-Petersson homeomorphisms on the real line. We denote  by $\WP_0(\Bbb C)$ the subset of all $h\in\WP(\Bbb C)$  with the normalized conditions (2.9), that is,
$$h(0)=0,\, h(\infty)=\infty,\, h(1)>0,\, \int_0^1|h'(t)|dt=1.\tag 8.1 $$
We also let $\Aff(\Bbb C)$ denote the set of all affine mappings $z\mapsto az+b$, $a\neq 0$.

\proclaim {Proposition 8.1}  The mapping $\Psi_1$ defined by $\Psi_1(h)=\log h'$ is a one-to-one map from $\Cal T_e\doteq\WP(\Bbb C)/\Aff(\Bbb C)$ into $H^{\frac 12}/\Bbb C$. The image $\Psi_1(\Cal T_e)$ is an open subset of $H^{\frac 12}/\Bbb C$.
\endproclaim
\demo{Proof} Clearly, $\Psi_1(h)=\log h'$ determines a one-to-one map $\Psi_1$ from $\Cal T_e$ into $H^{\frac 12}/\Bbb C$. We need to show that $\log h'_0$ is an interior point of $\Psi_1(\Cal T_e)$ for each $h_0\in\Cal T_e$. Let $w\in H^{\frac 12}/\Bbb C$ be given with small norm $\|w\|_{H^{\frac 12}}$. We need to find $h\in \Cal T_e$ with $\log h'=\log h'_0+w$.

Without loss of generality, we may assume that $h_0$ satisfies the normalized condition (8.1) so that $h_0\in\WP_0(\Bbb C)$. Consider the normalized Weil-Petersson curve $\Gamma_0=h_0(\hat\Bbb R)$ with the normalized arc-length parametrization $z=z_{\Gamma_0}$. Then there exists $b\in\hat\Cal T$ such that $z'=e^{ib}$. Consider the increasing homeomorphism $g_0$ on the real line determined by $z\circ g_0=h_0$ and set
 $$\tilde z(x)=\int_0^xe^{i(b(t)-i(w\circ g_0^{-1})(t))}dt.\tag 8.2$$
Since $\|w\|_{H^{\frac 12}}$ is small, and $g_0\in\WP_0(\Bbb R)$ since $\log g'_0=\log|h'_0|\in H^{\frac 12}$, we obtain  from Proposition 4.1 that $\|w\circ g_0^{-1}\|_{H^{\frac 12}}$ is also small. We conclude by the reasoning in section 6  that the equation (8.2) represents a Weil-Petersson  curve $\Gamma$. Set $h=\tilde z\circ g_0$ so that $h$ maps $\hat\Bbb R$ onto $\Gamma$. Then $h$ is locally absolutely continuous with
 $$h'=(\tilde z'\circ g_0)g'_0=(z'\circ g_0)g'_0e^w=h'_0e^w,$$
  which implies that $\log h'=\log h'_0+w$.  Consequently, $h\in\Cal T_e$ is the required mapping.
\quad$\square$
\enddemo

\proclaim{Proposition 8.2} There is a one-to-one map from $\Cal T_e$ onto $\WP_0(\Bbb R)\times \hat\Cal T$.
\endproclaim
\demo{Proof} From the proof of Proposition 8.1, each $h\in\WP(\Bbb C)$ induces a $g\in\WP(\Bbb R)$ and a $b\in\hat\Cal T$ such that $h=z\circ g$ maps $\hat\Bbb R$ onto a Weil-Petersson curve $\Gamma$ whose parametrization $z(s)$  by the arc-length $s\in\Bbb R$ satisfies $z'(s)=e^{ib(s)}$. This induces a one-to-one map $\Psi_2$ from $\Cal T_e$ onto $\WP_0(\Bbb R)\times \hat\Cal T$ by letting $\Psi_2(h)=(g, b)$. Actually, replacing $h$ by $\tilde h$ defined as
$$\tilde h(x)=\frac{|h(1)-h(0)|}{h(1)-h(0)}\frac{h(x)-h(0)}{\int_0^1|h'(t)|dt}$$
if necessary, we may assume that each $h\in\Cal T_e$ satisfies the normalized condition (8.1) so that $h\in\WP_0(\Bbb C)$. Then the corresponding function $g\in \WP(\Bbb R)$ satisfies the normalized condition $g(0)=0$, $g(1)$=1 so that $g\in\WP_0(\Bbb R)$.
\quad$\square$\enddemo

By means of Propositions 8.1 and 8.2, the generalized Weil-Petersson Teichm\"uller space $\Cal T_e$ can be endowed with two manifold structures. The following result says that they are topologically equivalent.

\proclaim{Proposition 8.3} The mapping $\hat\Psi\doteq\Psi_1\circ\Psi^{-1}_2$ is a homeomorphism from $\WP_0(\Bbb R)\times \hat\Cal T$ onto its image domain $\Psi_1(\Cal T_e)$ in $H^{\frac 12}/\Bbb C$.
\endproclaim
\demo{Proof} For $(g, b)\in\WP_0(\Bbb R)\times \hat\Cal T$, $\hat\Psi(g, b)=\log h'$, where $h=z\circ g$ with $z=z_b$ is the normalized arc-length parametrization of the normalized Weil-Petersson curve $\Gamma_b$. Thus,
$$\log h'=\log (z'\circ g)+\log g'=ib\circ g+\log g'.\tag 8.3$$
We conclude by Proposition 7.1 that $\hat\Psi$ is continuous from  $\WP_0(\Bbb R)\times \hat\Cal T$ into $H^{\frac 12}/\Bbb C$. Conversely, from (8.3) we obtain $\log g'=\Re\log h'$, and $b=\Im(\log h')\circ g^{-1}$, which implies by Proposition 7.1 again that $\hat\Psi^{-1}$ is continuous on $\Psi_1(\Cal T_e)$.
\quad$\square$\enddemo

\noindent
{\bf Proof of Theorem 2.5.}\quad  The correspondence $u\mapsto iu$ induces a homeomorphism from $\hat\Cal T_e$ onto the open subset $\Psi_1(\Cal T_e)$ in $H^{\frac 12}/\Bbb C$, which is contractible by Proposition 8.3. Consequently,  $\hat\Cal T_e$ is a contractible domain in $H^{\frac 12}/\Bbb C$. \quad$\square$

\vskip 0.2 cm
\noindent {\bf Acknowledgements} \quad
 The authors would like to thank the  referee for a  careful reading of the manuscript and for several corrections which  improves the presentation of the paper.


\enddocument